\documentclass{revtex4-2}
\usepackage{amsmath,amssymb}
\usepackage{graphicx}
\usepackage{fancyvrb}

\newtheorem{theorem}{Theorem}

\begin{document}
\title{Rigorous lower bounds for the domains of definition of extended generating
functions}
\author{B\'{e}la Erd\'{e}lyi, Jens Hoefkens, and Martin Berz}
\address{Department of Physics and Astronomy and National Superconducting Cyclotron\\
Laboratory, Michigan State University, East Lansing, MI 48824}
\date{\today}

\begin{abstract}
The recently developed theory of extended generating functions of symplectic
maps are combined with methods to prove invertibility via high-order Taylor
model methods to obtain rigorous lower bounds for the domains of definition
of generating functions. The examples presented suggest that there are
considerable differences in the size of the domains of definition of the
various types of generating functions studied. Furthermore, the types of
generating functions for which, in general, large domains of definitions can
be guarantied were identified. The resulting domain sizes of these
generators are usually sufficiently large to be useful in practice, for
example for the long-term simulations of accelerator and other Hamiltonian
system dynamics.
\end{abstract}

\maketitle

\section{Introduction}

It is well known that coordinate transformations of Hamiltonian systems that
keep the Hamiltonian structure intact are called canonical transformations,
or, in more modern terminology, symplectic maps \cite
{AIEP108book,goldstein,sanzserna}. The time evolutions of Hamiltonian
systems are also symplectic. Any symplectic map can be represented in terms
of a single scalar function, called the generating function. Until recently,
only the four conventional Goldstein-type generators were well known \cite
{goldstein}. Following the introduction of extended generating functions, it
has been shown that to each symplectic map, infinitely many generating
function types can be constructed \cite{prstabsymplectification}. In certain
applications, such as long term simulation of accelerators and other
Hamiltonian systems, it is important to maintain symplecticity during
tracking \cite{gjaja,abell,prstabsymplectification}. One available method to
achieve this is the generating function based symplectic tracking \cite
{caprimap,prstabsymplectification,warnock}.

In principle, any of the valid generators could be used for tracking of a
given Hamiltonian system. Although it can be shown that for any globally
defined symplectic map, global generators can always be found \cite{mythesis}%
, in practice their construction is not straightforward. Indeed, as shown in 
\cite{AIEP108book,prstabsymplectification}, the representation of symplectic
maps through generators commonly available in practice are often only
locally valid. Frequently, the purpose of the simulations is to estimate the
region of space where stable particle orbits can be found, the so-called
non-wandering set, or dynamic aperture in the beam dynamics terminology.
Hence, a sizable phase space region must be covered by tracking, and if the
generating function method is used, the generating function must be defined
at least in that region. However, sofar nothing has been known about the
size of the domains of definition of generating functions. The necessity to
study this interesting problem has been recognized in \cite
{caprimap,warnock,chansus,prstabsymplectification}.

The purpose of this paper is to show that, using sophisticated high-order
Taylor model based methods \cite{rdasf,ACAInv}, it is possible to rigorously
prove lower bounds on the size of the domains of definition of any type of
generating function. The resulting domains often enclose the dynamic
apertures of the examples studied.

The structuring of the paper is as follows: in Section \ref{sect:genfunc}
the extended generating function theory is reviewed, and in Section \ref
{sect:numtools} a brief description of the verified numerical tools utilized
for bounding the domains of definition is presented. Finally, the
performance of the methods is illustrated by four examples in Section \ref
{sect:ex}: two polynomial maps in two and four dimensions, respectively, the
Fermi-Pasta-Ulam system \cite{fpu}, and a representative example taken from
beam physics.

\section{The extended generating function theory}

\label{sect:genfunc}

Following the exposition of \cite{prstabsymplectification}, we regard every
map as a column vector. Let 
\begin{equation}
\alpha =\left( 
\begin{tabular}{c}
$\alpha _{1}$ \\ 
$\alpha _{2}$%
\end{tabular}
\right)
\end{equation}
be a diffeomorphism of a subset of ${\Bbb R}^{4n}$ onto its image. Notice
that $\alpha _{i},$ $i=1,2$, are the first $2n$ and second $2n$ components
of $\alpha $. This entails that $\alpha _{i}:U\subset {\Bbb R}%
^{4n}\rightarrow V_{i}\subset {\Bbb R}^{2n}$. Let 
\begin{equation}
%TCIMACRO{\func{Jac}}%
%BeginExpansion
\mathop{\rm Jac}%
%EndExpansion
\left( \alpha \right) =\left( 
\begin{array}{cc}
A & B \\ 
C & D
\end{array}
\right)
\end{equation}
be the $4n\times 4n$ Jacobian of $\alpha $, split into $2n\times 2n$ blocks.
Let 
\begin{equation}
\tilde{J}_{4n}=\left( 
\begin{array}{cc}
J_{2n} & 0_{2n} \\ 
0_{2n} & -J_{2n}
\end{array}
\right) ,
\end{equation}
where 
\begin{equation}
J_{2n}=\left( 
\begin{array}{cc}
0_{n} & I_{n} \\ 
-I_{n} & 0_{n}
\end{array}
\right)
\end{equation}
with the unit matrix $I_{n}$ of appropriate dimension. A map $\alpha $ is
called conformal symplectic if 
\begin{equation}
\left( 
%TCIMACRO{\func{Jac}}%
%BeginExpansion
\mathop{\rm Jac}%
%EndExpansion
\left( \alpha \right) \right) ^{T}J_{4n}%
%TCIMACRO{\func{Jac}}%
%BeginExpansion
\mathop{\rm Jac}%
%EndExpansion
\left( \alpha \right) =\mu \tilde{J}_{4n},  \label{eq:confsymdef}
\end{equation}
where $\mu $ is a non-zero real constant \cite{banyagabook}. Also, we denote
by ${\cal I}$ the identity map of appropriate dimension. A map ${\cal M}$ is
called symplectic if its Jacobian $M$ satisfies the symplectic condition 
\cite{dragtlectures}, that is 
\begin{equation}
M^{T}JM=J.
\end{equation}
Henceforth, we will consider symplectic maps that are origin preserving,
i.e. maps around a fixed point. We call a map a gradient map if it has
symmetric Jacobian $N$. It is well known that, at least over simply
connected subsets, gradient maps can be written as the gradient of a
function (hence the name) \cite{AIEP108book}, that is 
\begin{equation}
N=%
%TCIMACRO{\limfunc{Jac}}%
%BeginExpansion
\mathop{\rm Jac}%
%EndExpansion
\left( \nabla F\right) ^{T}.
\end{equation}
($\nabla F$ is regarded as a row vector \cite{marsdenbook2}.) The function $%
F $ is called the potential of the map.

The best way to formulate the main result of the extended generating
function theory is the following theorem \cite{prstabsymplectification}:

\begin{theorem}[Existence of extended generating functions]
\label{thm:extgen} Let ${\cal M}$ be a symplectic map. Then, for every point 
$z$ there is a neighborhood of $z$ such that ${\cal M}$ can be represented
by functions $F$ via the relation 
\begin{equation}
\left( \nabla F\right) ^{T}=\left( \alpha _{1}\circ \left( 
\begin{tabular}{c}
${\cal M}$ \\ 
${\cal I}$%
\end{tabular}
\right) \right) \circ \left( \alpha _{2}\circ \left( 
\begin{tabular}{c}
${\cal M}$ \\ 
${\cal I}$%
\end{tabular}
\right) \right) ^{-1},  \label{eq:fdef}
\end{equation}
where $\alpha =\left( \alpha _{1},\alpha _{2}\right) ^{T}$ is any conformal
symplectic map such that 
\begin{equation}
\det \left( C\left( {\cal M}\left( z\right) ,z\right) \cdot Mz+D\left( {\cal %
M}\left( z\right) ,z\right) \right) \neq 0.
\end{equation}
Conversely, let $F$ be a twice continuously differentiable function. Then,
the map ${\cal M}$ defined by 
\begin{equation}
M=\left( NC-A\right) ^{-1}\left( B-ND\right)
\end{equation}
is symplectic. The matrices $A,B,C,D,M,$ and $N$ are defined above.
\end{theorem}

The function $F$ is called the generating function of type $\alpha $ of $%
{\cal M}$, and denoted by $F_{\alpha ,{\cal M}}$. Theorem \ref{thm:extgen}
says that, once the generator type $\alpha $ is fixed, locally there is a
one-to-one correspondence between symplectic maps and scalar functions,
which are unique up to an additive constant. The constant can be normalized
to zero without loss of generality. Due to the fact that there exist
uncountably many maps of the form (\ref{eq:confsymdef}), to each symplectic
map one can construct infinitely many generating function types.

We note that (\ref{eq:fdef}) cannot be simplified, due to the fact that the
entries in the equations have different dimensions. For instance, $\alpha
_{i}$ maps ${\Bbb R}^{4n}\rightarrow {\Bbb R}^{2n}$ and $\left( 
\begin{tabular}{cc}
${\cal M}$ & ${\cal I}$%
\end{tabular}
\right) ^{T}$ maps ${\Bbb R}^{2n}\rightarrow {\Bbb R}^{4n}$. For more
details see \cite{prstabsymplectification}.

Therefore, for a given ${\cal M}$, the domain of definition of the generator
of type $\alpha $ is the domain of invertibility of 
\begin{equation}
\left( \alpha _{2}\circ \left( 
\begin{tabular}{c}
${\cal M}$ \\ 
${\cal I}$%
\end{tabular}
\right) \right) ,  \label{eq:inv}
\end{equation}
according to (\ref{eq:fdef}). Thus, finding the domain of definition of a \
generator is equivalent to finding the domain of invertibility of (\ref
{eq:inv}). A large class of generator types, which are utilized in practice,
is the generators associated with linear maps $\alpha $. These can be
organized into equivalence classes $\left[ S\right] $, associated with 
\begin{equation}
\alpha =\left( 
\begin{array}{cc}
-JL^{-1} & J \\ 
\frac{1}{2}\left( I+JS\right) L^{-1} & \frac{1}{2}\left( I-JS\right)
\end{array}
\right) ,  \label{eq:implementalpha}
\end{equation}
and represented by arbitrary symmetric matrices $S$ \cite
{prstabsymplectification}. The matrix $L$ stands for the linear part of $%
{\cal M}$.

In the next section we present the methods utilized to prove lower bounds of
the domains of definition of a variety of generators, or equivalently
symmetric matrices $S$, for several given ${\cal M}$s.

\section{Rigorous Numerical Analysis with Taylor Models}

\label{sect:numtools}

The relative accuracy of floating point number representations on computers
is limited, usually to around 16 significant digits. However, interval
analysis takes this into account and allows the computation of guaranteed
enclosures for computational results: the result of every interval
computation consists of two numbers. Thus, if $x_{r}$ is the analytical
result of some computation $C$, and $R=[x_{l},x_{u}]$ is the corresponding
result of interval analysis, then $x_{r}\in R$. On the other hand, interval
analysis can also be seen as a computerization of set theory: starting with $%
I=[x_{l},x_{u}]$ and some computation $C$, the result $C(I)=[y_{l},y_{u}]$
contains the result of the computation $C$ applied to {\em every} number
contained in $I$. Further details on interval analysis can be found in \cite
{moore66,Hansen69,moore79}.

Recently, Taylor models have been developed \cite{rdasf} as a combination of
the well known map approach \cite{AIEP108book} and interval analysis to
obtain guaranteed enclosures of functional descriptions. Details on
applications of these methods in beam physics are given in \cite
{rdapac97,VICAP00}. Loosely speaking, a Taylor model is a set of maps that
are, in some yet-to-be specified way, close to a reference polynomial. We
consider a box ${\mbox{\boldmath$D$}}\subset {\Bbb R}^{v}$ containing $x_{0}$
and assume $P_{n}:{\mbox{\boldmath$D$}}\rightarrow {\Bbb R}^{w}$ to be a
polynomial of order $n$, where $n,v,w$ are all integers. If we then take an
open non-empty set $R\subset {\Bbb R}^{w}$, we call the collection $\left(
P_{n},x_{0},{\mbox{\boldmath$D$}},R\right) $ a {\em Taylor model of order $n$
with reference point $x_{0}$ over ${\mbox{\boldmath$D$}}$}.

We say that $f$ is contained in a Taylor model $T=\left( P_{n},x_{0},{%
\mbox{\boldmath$D$}},R\right) $ if $P_{n}(x)-f(x)\in R$ for all $x\in {%
\mbox{\boldmath$D$}}$ and the $n$-th order Taylor series of $f$ around $%
x_{0} $ equals $P_{n}$.

Methods have been developed that allow mathematical operations on Taylor
models that preserve the inclusion relationships. For example, for two given
Taylor models, $T_{1}$ and $T_{2}$, methods have been developed that compute
Taylor models for the sum $S$, and the product $P$, of $T_{1}$ and $T_{2}$,
i.e. 
\begin{eqnarray*}
f_{1}\in T_{1},f_{2}\in T_{2} &\Rightarrow &(f_{1}+f_{2})\in S \\
f_{1}\in T_{1},f_{2}\in T_{2} &\Rightarrow &(f_{1}\cdot f_{2})\in P.
\end{eqnarray*}
Similarly, the other elementary operations and intrinsic functions have been
extended to Taylor models such that the fundamental inclusion properties of
Taylor models are maintained. For the exponential function, this means that
we have developed algorithms that compute a Taylor model $T_{E}$ from a
given Taylor model $T$ such that 
\begin{eqnarray*}
f\in T &\Rightarrow& \exp\left(f\right) \in T_E
\end{eqnarray*}

Lastly, the antiderivation, which is essentially the integration operation,
extends naturally to Taylor models. It allows us to compute a Taylor model $%
T_{I}$ from a given Taylor model $T$ such that $T_{I}$ contains primitives
for each function $f$ contained in $T$. Since the antiderivation does not
fundamentally differ from other intrinsic functions on the set of Taylor
models, it is often used in fundamental Taylor model algorithms. An
important application of the antiderivation will be presented in (\ref
{eq:picard}).

More information on arithmetic and intrinsic functions on Taylor models can
be found in \cite{rcdep,rdasf,rdaic}.

One major advantage of Taylor model methods over regular interval analysis
is that the accuracy of the enclosures obtained by Taylor models scales with
the $(n+1)$-st order of the domain size \cite{RCInv}. Thus, the Taylor model
approach is of particular advantage in the combination of Taylor models with
high-order map codes like {\sc cosy infinity} \cite{cosyman8} and gives
small guaranteed enclosures even for complicated high-order maps in many
variables.

\subsection{Guaranteed Invertibility of Maps}

The Inverse Function Theorem makes local statements about invertibility of
maps at points $x_{0}$. On the other hand, interval methods allow us to {\em %
test} whether a given map is invertible over a whole region: assuming that
we have a criterion that has to be satisfied for {\em each point} in a
region for the function to be invertible, we can use interval methods to
test the whole region {\em at once}.

One such criterion is given below in Theorem \ref{thm:invertibility}. While
it does not lend itself very well to pure interval arithmetic, since normal
interval methods model the inherent structure of the matrix $M$ rather
inefficiently, it has been shown in \cite{ACAInv} that the Taylor model
approach can be used to make full use of the specific row-structure of $M$.

\begin{theorem}[Invertibility from First Derivatives]
\label{thm:invertibility} Let $\mbox{\boldmath $B$}\subset {\Bbb R}^{v}$ be
a box and $f:\mbox{\boldmath $B$}\rightarrow {\Bbb R}^{v}$ a ${\cal C}^{1}$
function. Assume that the matrix 
\begin{equation}
M=\left( 
\begin{array}{ccc}
\frac{\partial f_{1}}{\partial x_{1}}(\chi _{1}) & \cdots & \frac{\partial
f_{1}}{\partial x_{v}}(\chi _{1}) \\ 
\vdots &  & \vdots \\ 
\frac{\partial f_{v}}{\partial x_{1}}(\chi _{v}) & \cdots & \frac{\partial
f_{v}}{\partial x_{v}}(\chi _{v})
\end{array}
\right)
\end{equation}
is invertible for every choice of $\chi _{1},\ldots ,\chi _{v}\in 
\mbox{\boldmath
$B$}$. Then $f$ has a ${\cal C}^{1}$-inverse defined on the whole set $f(%
\mbox{\boldmath $B$})$, where $f(\mbox{\boldmath $B$})$ denotes the range of 
$f$ over $\mbox{\boldmath $B$}$.
\end{theorem}

It has been shown in \cite{ACAInv} that a Taylor model implementation of
Theorem \ref{thm:invertibility} can successfully determine invertibility of
complicated high-dimensional functions over large regions $%
\mbox{\boldmath
$B$}$. We do, however, point out that an application of this theorem cannot
disprove invertibility of a given map.

\subsection{Verified Enclosures of Flows}

Another important application of Taylor models stems from their
applicability to numerical ODE integration. To illustrate how this works, we
consider the initial value problem 
\begin{equation}
x^{\prime }=f(t,x)\;\;\mbox{ and }\;\;x(t_{0})=x_{0}.  \label{eq:ivp}
\end{equation}
It is a well known fact that the solution to (\ref{eq:ivp}) can be obtained
as the fixed point of the Picard operator ${\cal O}$, defined by 
\begin{equation}
{\cal O}\left( x\right) =x_{0}+\int_{t_{0}}^{t}f\left( \tau ,x\right) d\tau .
\label{eq:picard}
\end{equation}
Using the previously mentioned antiderivation, the operator ${\cal O}$ can
be extended to Taylor models and yields an algorithm that allows the
computation of verified {\em enclosures of flows} of (\ref{eq:ivp}). This
approach has been presented in \cite{rdaint,makinophd} and has recently been
used successfully in a variety of applications ranging from solar system
dynamics \cite{ASNL00} to beam physics. Unlike methods that use regular
interval computations to enclose the final conditions, the Taylor model
approach avoids the wrapping effect to very high order. It is therefore
capable of propagating extended initial regions over large integration
intervals.

If we denote the solution of (\ref{eq:ivp}) by ${\cal M}\left(
t,t_{0},x_{0}\right) $, we consider the corresponding matrix initial value
problem 
\begin{equation}
Y^{\prime }=\left( \frac{\partial f}{\partial x}\left( t,{\cal M}\left(
t,t_{0},x_{0}\right) \right) \right) \cdot Y,\;\;Y(t_{0})=I.
\label{eq:adjoint}
\end{equation}
We use Taylor model methods to compute Taylor models containing the map $%
{\cal M}(t_{f},t_{0},x_{0})$ at the final time $t_{f}$. Moreover, since $%
\frac{\partial {\cal M}}{\partial x_{0}}$ satisfies (\ref{eq:adjoint}), we
also obtain Taylor models containing the derivatives $\frac{\partial {\cal M}%
}{\partial x_{0}}(t_{f},t_{0},x_{0})$ at that final point. The latter is
required for a determination of invertibility of (\ref{eq:inv}) based on
Theorem \ref{thm:invertibility}.

\section{Examples}

\label{sect:ex}

\subsection{A cubic $2D$ symplectic map}

As a first example we consider a two dimensional cubic polynomial, which is
exactly symplectic. It has the following form: 
\begin{equation}
{\cal M}={\cal N}\circ {\cal L}\text{,}  \label{eq:quadsymp}
\end{equation}
where 
\begin{equation}
{\cal L}=\left( 
\begin{array}{cc}
\cos \frac{\pi }{3} & \sin \frac{\pi }{3} \\ 
-\sin \frac{\pi }{3} & \cos \frac{\pi }{3}
\end{array}
\right) ,
\end{equation}
and 
\begin{equation}
{\cal N}\left( 
\begin{tabular}{c}
$q$ \\ 
$p$%
\end{tabular}
\right) =\left( 
\begin{tabular}{c}
$q-3\left( q+p\right) ^{3}$ \\ 
$p+3\left( q+p\right) ^{3}$%
\end{tabular}
\right) .
\end{equation}
This simple example even lends itself to an analytical treatment. Applying
Theorem \ref{thm:invertibility}, we obtain the following determinant in the
four variables $q_{1},p_{1},q_{2},p_{2}$: 
\begin{eqnarray}
&&1+9\frac{\left[ \left( 1-\sqrt{3}\right) q_{1}+\left( 1+\sqrt{3}\right)
p_{1}\right] ^{2}}{8\left( 1+\sqrt{3}\right) ^{2}}\left[ \left( 2+\sqrt{3}%
\right) \left( s_{12}+1\right) +s_{22}\right]  \nonumber \\
&&+9\frac{\left[ \left( 1-\sqrt{3}\right) q_{1}+\left( 1+\sqrt{3}\right)
p_{1}\right] ^{2}}{8\left( 1+\sqrt{3}\right) ^{2}}\left[ \left( 2+\sqrt{3}%
\right) \left( s_{12}-1\right) +\left( 7+4\sqrt{3}\right) s_{11}\right] ,
\label{eq:det}
\end{eqnarray}
where 
\begin{equation}
S=\left( 
\begin{array}{cc}
s_{11} & s_{12} \\ 
s_{12} & s_{22}
\end{array}
\right) .
\end{equation}
Theorem \ref{thm:invertibility}\ guaranties globally defined inverses if (%
\ref{eq:det}) is positive definite. Clearly, this is the case if 
\begin{eqnarray}
s_{22} &\leq &\left( 2+\sqrt{3}\right) \left( s_{12}+1\right) \\
s_{11} &\leq &\frac{2+\sqrt{3}}{7+4\sqrt{3}}\left( s_{12}-1\right) ,
\end{eqnarray}
for any $s_{12}$. Therefore, we are able to prove that there is a large
class of generators that is globally defined for this symplectic map.

\begin{figure}
\centering
\includegraphics[scale=0.65]{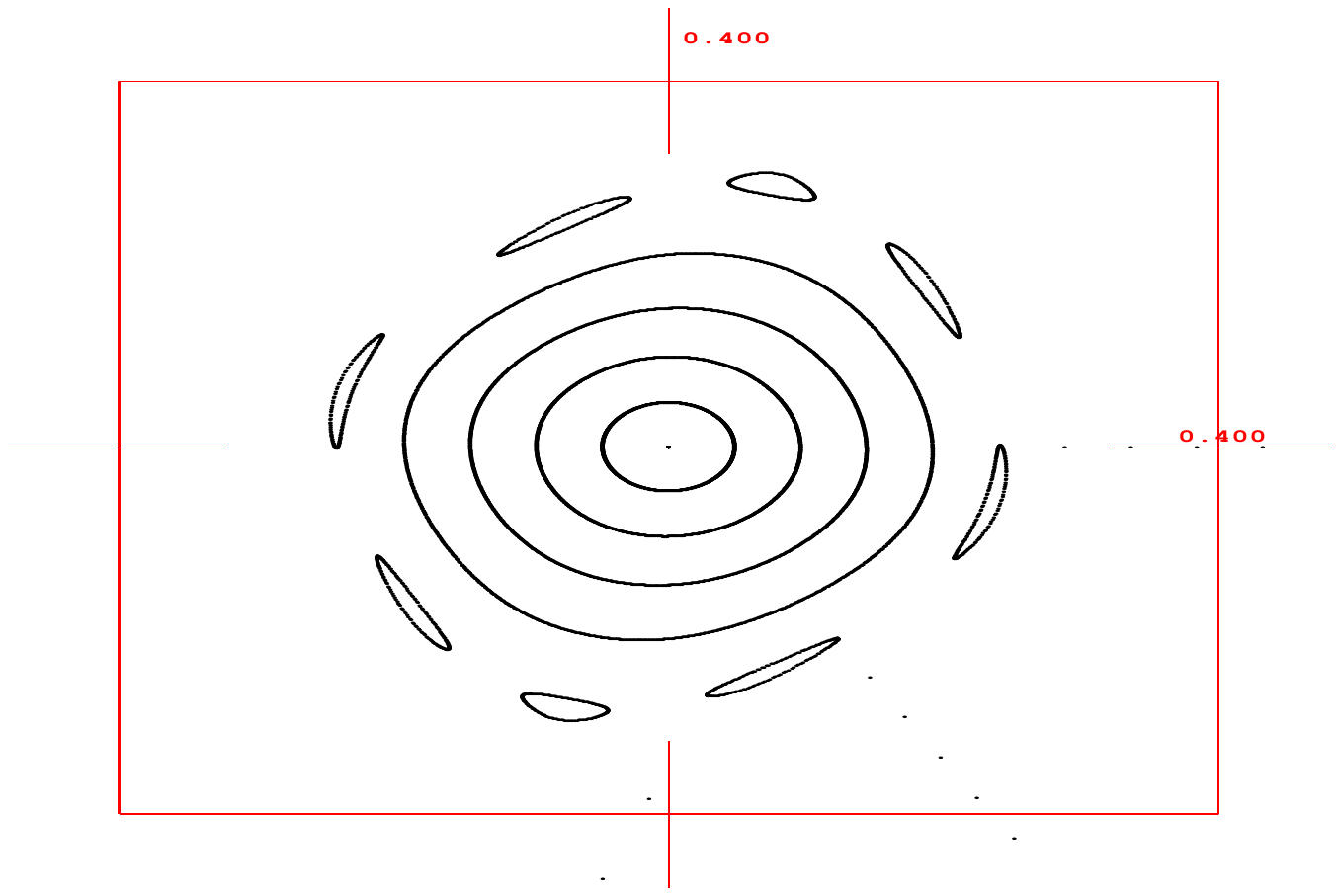}
\caption
{Tracking picture of the cubic two dimensional symplectic map, and  the box of 
guarantied invertibility of the generator associated with $S=0$.}
\label{fig:2dpoly}
\end{figure}

However, these solutions are not unique, and we are interested mainly in
proving invertibility in finite regions that enclose the dynamic aperture.
For these cases, the Taylor model based method described in Section \ref
{sect:numtools} is well suited. For example, employing the Taylor method for
the $S=0$ case, the tracking picture shown in Fig. \ref{fig:2dpoly} was
obtained. The box surrounding the tracking is the region for which
invertibility of the corresponding generator can be proved. Clearly, it goes
well beyond the dynamic aperture, which is approximately $0.2$. To check
invertibility for other generators, we created a few random generator types,
based on random symmetric matrices with entries uniformly distributed in $%
\left[ -10,10\right] $. The results are presented in Table \ref{table:2dpoly}%
.

\subsection{A $4D$ symplectic map}

The second example is somewhat more complex. While staying in the category
of polynomial symplectic maps, the dimensions were doubled $\left( 4D\right) 
$, and to ensure exact symplecticity, we composed two types of maps that are
known to be symplectic, namely rotations 
\begin{equation}
{\cal R}=\left( 
\begin{array}{cc}
R_{1} & 0 \\ 
0 & R_{2}
\end{array}
\right) ,
\end{equation}
where 
\begin{equation}
R_{i}=\left( 
\begin{array}{cc}
\cos \frac{\pi }{3} & \sin \frac{\pi }{3} \\ 
-\sin \frac{\pi }{3} & \cos \frac{\pi }{3}
\end{array}
\right) ,
\end{equation}
$i=1,2$, and kicks 
\begin{equation}
{\cal K}_{i}\left( 
\begin{tabular}{c}
$\vec{q}$ \\ 
$\vec{p}$%
\end{tabular}
\right) =\left( 
\begin{tabular}{c}
$q_{i}$ \\ 
$p_{i}+\frac{\partial f}{\partial q_{i}}\left( \vec{q}\right) $%
\end{tabular}
\right) ,
\end{equation}
where $\vec{p}=\left( p_{1},p_{2}\right) ^{T}$, and $f$ is any polynomial of 
$\vec{q}=\left( q_{1},q_{2}\right) ^{T}$.

\begin{figure}
\centering
\includegraphics[scale=0.65]{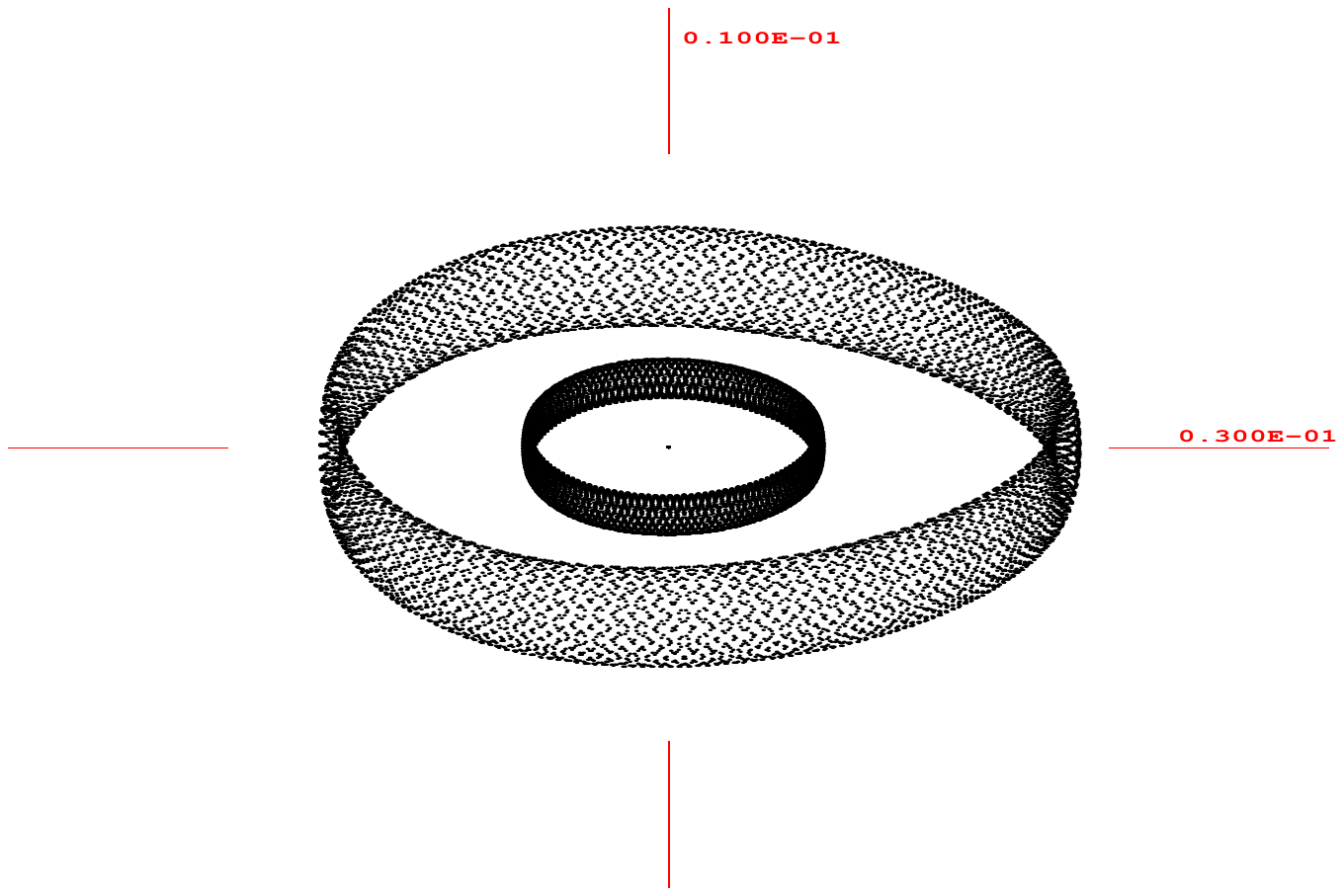}
\caption
{$(q_1,p_1)$ tracking picture of the four dimensional symplectic polynomial, 
for two initial conditions launched along the $q_1$ axis close to the origin.}
\label{fig:4dpoly-1}
\end{figure}

\begin{figure}
\centering
\includegraphics[scale=0.3]{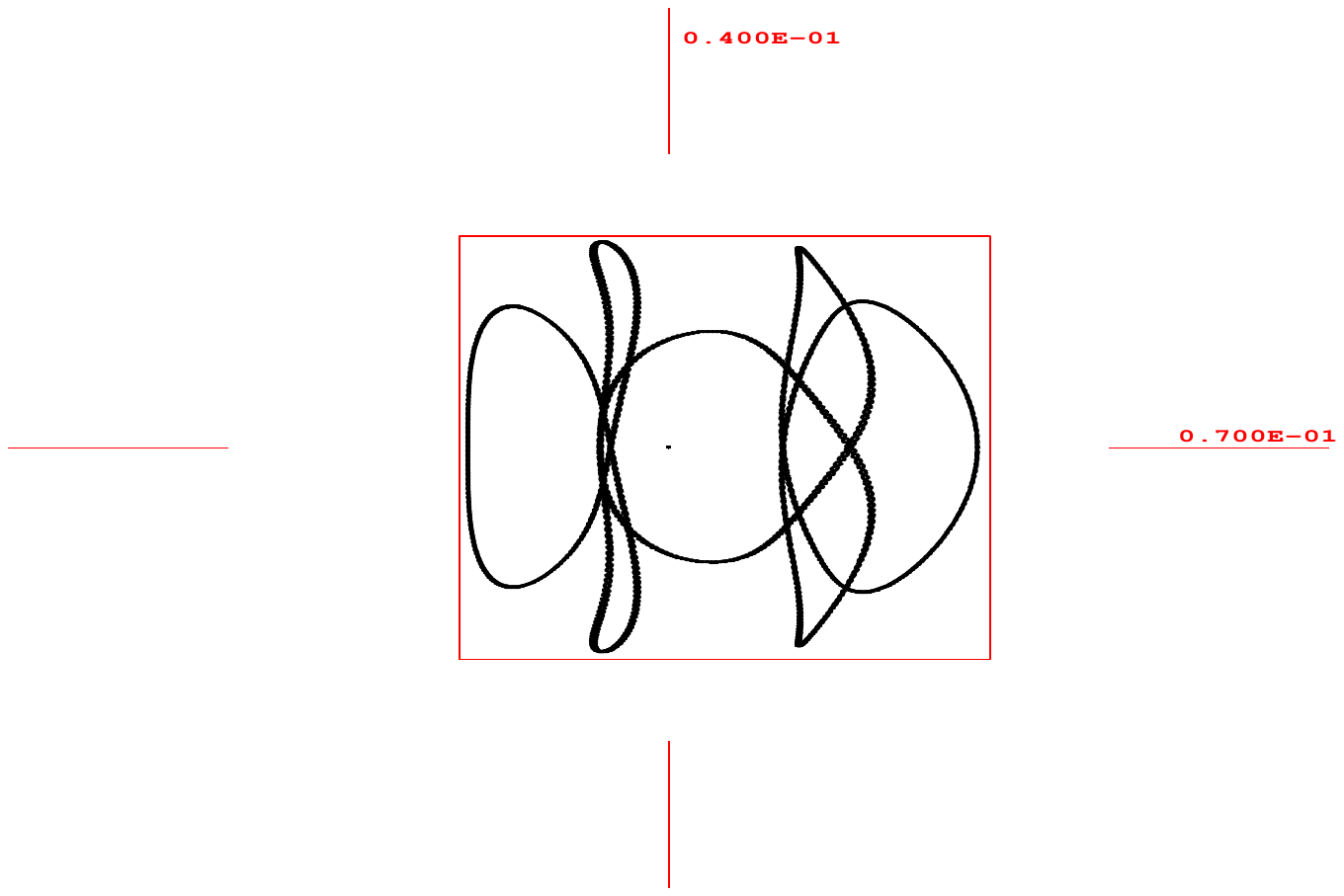}
\includegraphics[scale=0.3]{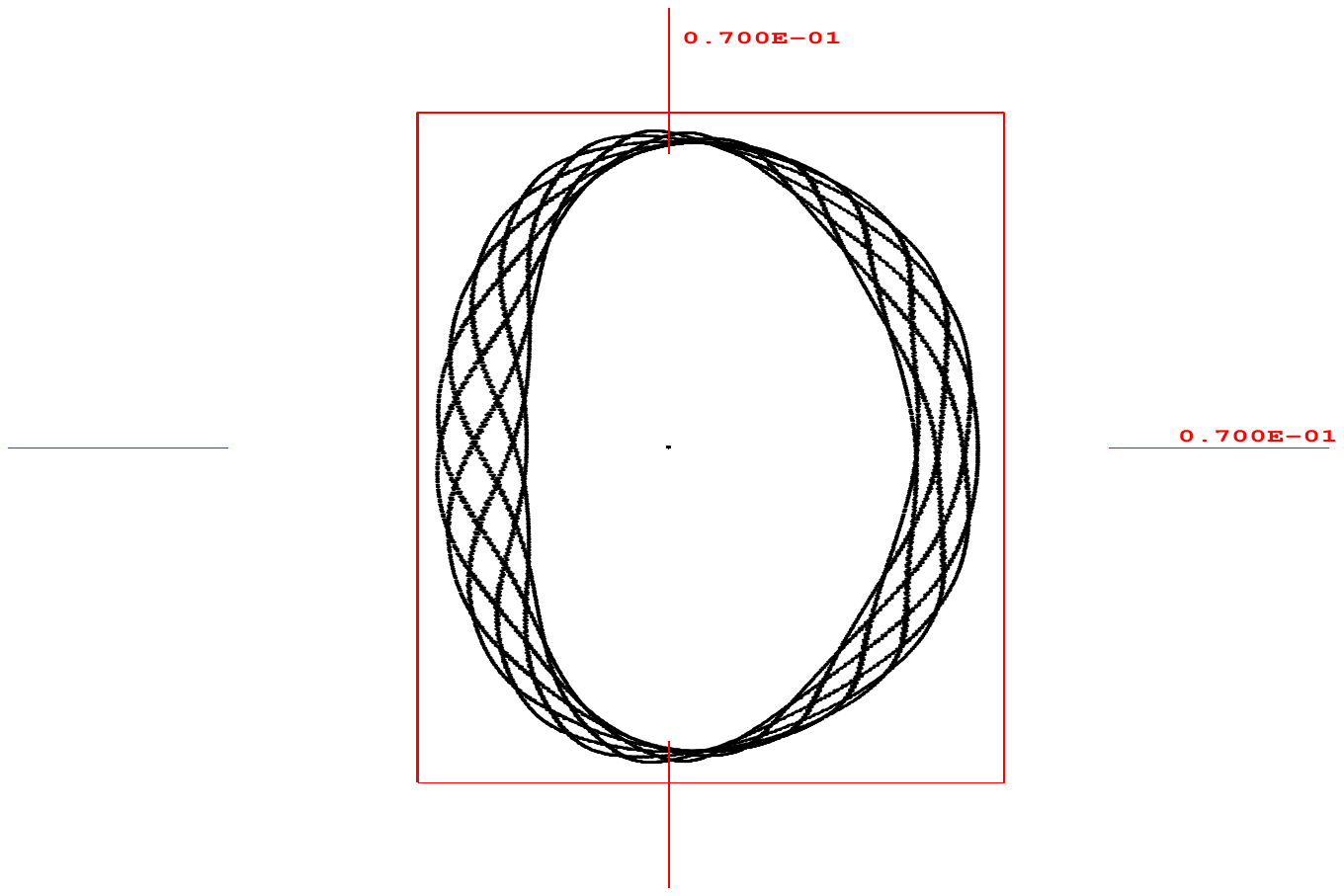}
\caption
{$(a)$ $(q_1,p_1)$, and $(b)$ $(q_2,p_2)$ tracking pictures of the four dimensional 
symplectic polynomial for two particles (launched along the $q_1$ and $q_2$ axes 
respectively) close to the dynamic aperture, and  the box of 
guarantied invertibility of the generator associated with $S=0$.}
\label{fig:4dpoly-2}
\end{figure}

For a typical example (with $f\left( q_{1},q_{2}\right)
=1.5q_{1}q_{2}+q_{1}^{2}+2\left( q_{1}+q_{2}\right) ^{3}+1.45\left(
1.1q_{1}-0.5q_{2}\right) ^{4}$), we obtained the tracking pictures depicted
in Figs. \ref{fig:4dpoly-1} and \ref{fig:4dpoly-2}. Fig. \ref{fig:4dpoly-1}
shows some particles launched relatively close to the origin, where the
system is still weakly nonlinear, and the behavior of the particles is
similar to those in a real accelerator. However, the dynamic aperture is
much bigger, and lies somewhere in the more nonlinear region, close to the
particles tracked in Fig. \ref{fig:4dpoly-2}. As in the previous example, in
Fig. \ref{fig:4dpoly-2} the lower bounds for invertibility of the case $S=0$
are also shown. In fact, long term tracking shows that the dynamic aperture
is approximately $0.036$ in $q_{1}$ and $0.052$ in the $q_{2}$ directions.
Applying the Taylor method to a variety of generators we found the lower
bounds for the domains shown in Table \ref{table:4dpoly}. Notice that
several generators can be guaranteed to be defined beyond the dynamic
aperture.

\subsection{The Fermi-Pasta-Ulam system}

The first two examples have been fabricated such that the corresponding
symplectic maps are polynomial and known to arbitrary precision. This is not
usually the case for problems encountered in practice. The information that
is supposed to be known exactly are the Hamiltonian functions, or
equivalently, the corresponding ordinary differential equations. Since
closed form solutions are almost never available, the solutions are usually
approximated, for example by Taylor series. Using verified integration
methods, we can enclose the true solution by a Taylor model. The Taylor
model contains many symplectic (and non-symplectic) maps, one of which is
the exact solution. Therefore, applying the invertibility tests to the
Taylor model containing the true solution of the ODE, we can again prove
lower bounds for domains of generating functions.

\begin{figure}
\centering
\includegraphics[scale=1]{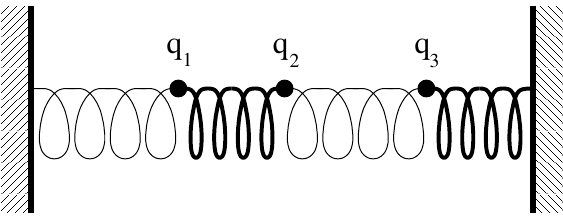}
\caption
{The Fermi-Pasta-Ulam system used in our example.}
\label{fig:springs}
\end{figure}

\begin{figure}
\centering
\includegraphics[scale=0.3]{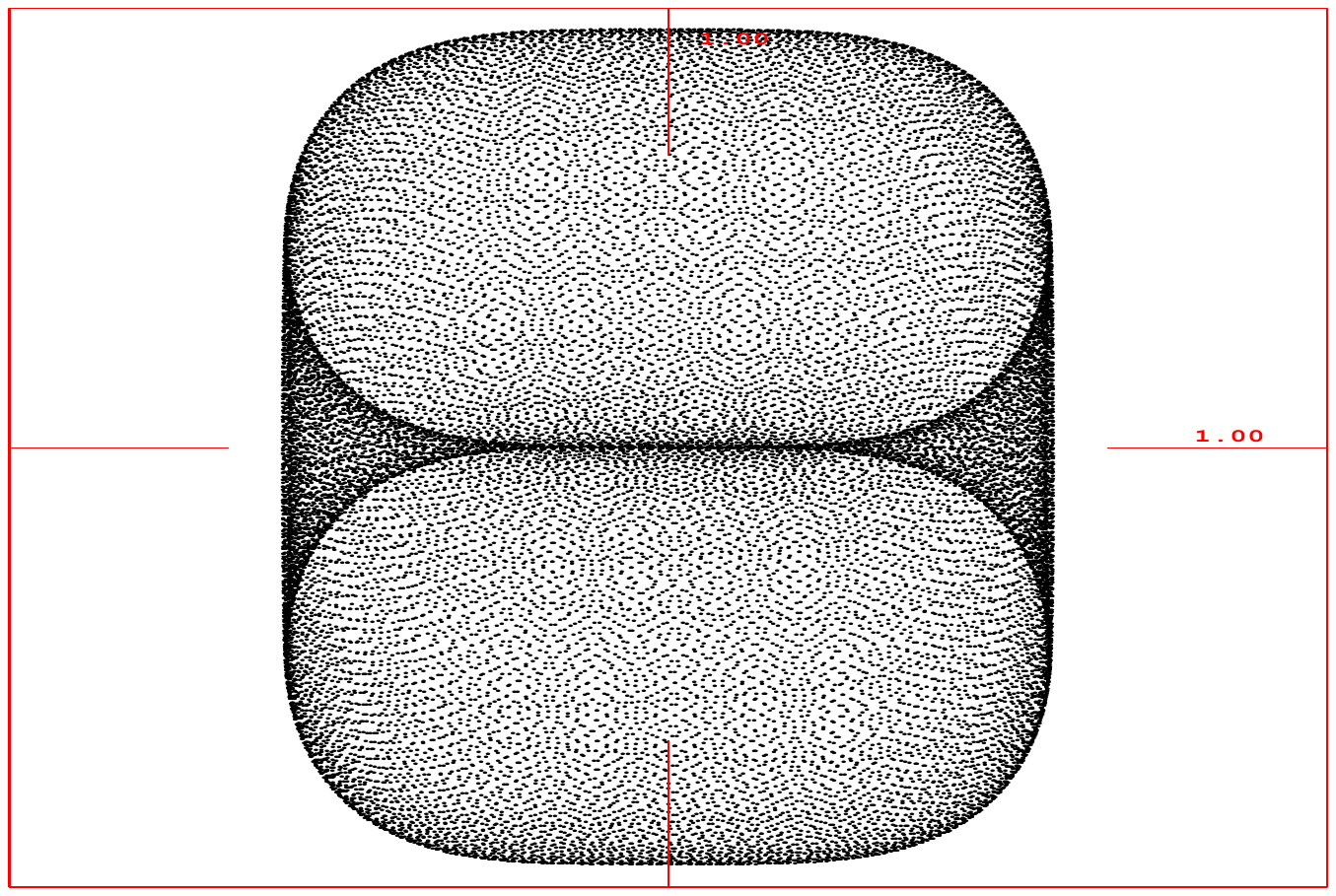}
\includegraphics[scale=0.3]{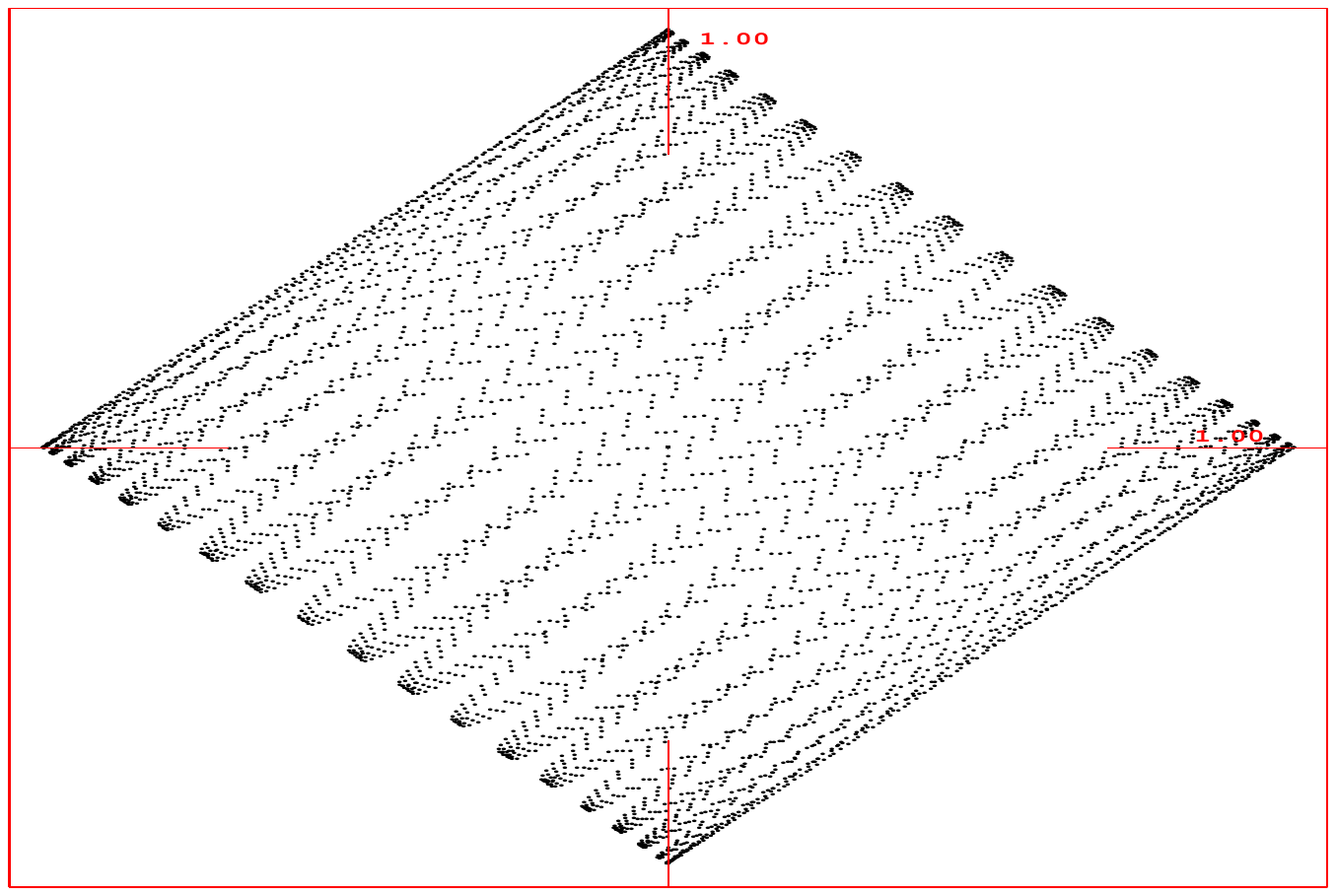}
\caption
{$(a)$ $(q_1,p_1)$, and $(b)$ $(q_1,q_2)$ tracking pictures of the half period map of 
the Fermi-Pasta-Ulam system for a particle launched along the $q_1$ axis. The box of 
guarantied invertibility of the generator associated with $S=0$ extends to at least 
$[-1,1]$ in every direction.}
\label{fig:fpu}
\end{figure}

A famous Hamiltonian system is the Fermi-Pasta-Ulam system shown in Fig. \ref
{fig:springs} \cite{fpu}. It is a paradigmatic model for nonlinear classical
systems, which is known to exhibit many interesting features, and stimulated
important developments in nonlinear dynamics. It is determined by the
Hamiltonian 
\begin{equation}
H\left( \vec{q},\vec{p}\right) =\frac{1}{2}\sum_{i=1}^{n}\left(
p_{2i-1}^{2}+p_{2i}^{2}\right) +\frac{k}{2}\sum_{i=1}^{n}\left(
q_{2i}-q_{2i-1}\right) ^{2}+\sum_{i=0}^{n}\left( q_{2i+1}-q_{2i}\right) ^{4}.
\label{eq:fpu-ham}
\end{equation}
For our specific example we took $n=2$ $\left( q_{4}=q_{5}=p_{4}=0\right) $,
and $k=5000$. Employing the Taylor model based verified integrator, we
obtained the Taylor model of one half period map of the flow of (\ref
{eq:fpu-ham}), and then tracked some particles by iteration of the map.
Typical phase space plots are presented in Fig. \ref{fig:fpu}. The domain of
the generator associated with $S=0$ extends to at least $\left[ -1,1\right] $
in every direction in phase space. Since there is no real dynamic aperture
in this case, to compare different generating function types a set of $1000$
random symmetric matrices were created with entries in $\left[ -1,1\right] $%
. Then the corresponding generator types were constructed. Furthermore, the
symmetric matrices were multiplied by constants between $\left[ 0,1\right] $
to investigate the percentage of generators associated to these matrices
that are still defined over $\left[ -1,1\right] $ in every direction in
phase space. As shown in Fig. \ref{fig:size-of-s}, there is a sharp increase
in the fraction of generators defined over the whole region of interest as
the ``size'' of the symmetric matrices decrease. If $S$ is small enough, the
generator type associated with it is defined over $\left[ -1,1\right] $.

\subsection{An accelerator cell}

\label{ssect:fodo}

\begin{figure}
\centering
\includegraphics[scale=0.65]{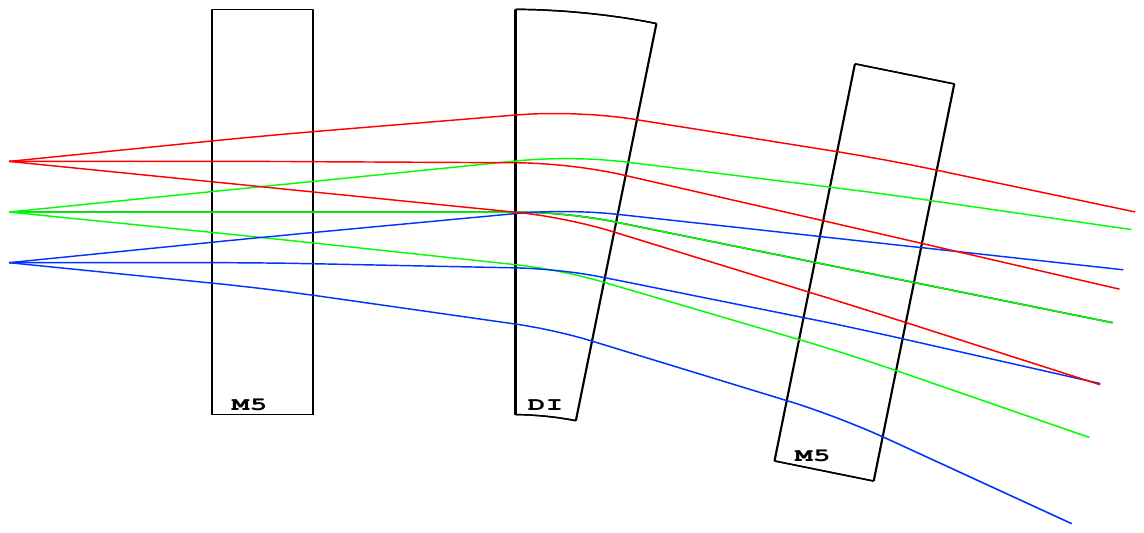}
\caption
{The accelerator cell used in the example of Subsection \ref{ssect:fodo}. The drifts are 
$1 m$ long, the quadrupoles and the dipole are $0.5 m$ long, and the dipole's radius 
of curvature is $2.5 m$. The first quadrupoles are defocusing in the $x$ direction, 
with strength $k=-0.0085$, and the reference particle is a $1 MeV$ proton.}
\label{fig:fodo}
\end{figure}

We conclude the section with an example of practical interest in accelerator
physics. The Hamiltonian representing an accelerator magnet, with the
arclength $s$ as independent variable, and on-energy particles is 
\begin{equation}
H\left( x,y,a,b\right) =-\left( 1+\frac{x}{\rho }\right) \sqrt{1-a^{2}-b^{2}}%
-\frac{q}{p_{0}}\left( 1+\frac{x}{\rho }\right) A_{s}\left( x,y\right) ,
\end{equation}
where $\rho $ is the curvature radius of the magnet, $A_{s}$ is the $s$
component of the vector potential, and $\left( \left( x,a);(y,b\right)
\right) $ are two pairs of canonically conjugate variables. For the sake of
computational simplicity, we assume that the magnetic fields are piecewise
independent of the arclength, in particular we use the sharp cutoff
approximation for fringe fields, as is commonly done in beam physics.
Specifically, for field free space $\rho =\infty $ and $A_{s}=0$, for a
homogeneous dipole magnet 
\begin{equation}
A_{s}=-\frac{B_{0}}{2}\left( x+\rho \right) ,
\end{equation}
where $qB_{0}/p_{0}=1/\rho $, for a quadrupole $\rho =\infty $ and 
\begin{equation}
A_{s}=-\frac{k}{2}\left( x^{2}-y^{2}\right) ,
\end{equation}
where $k$ is the quadruple strength ($k>0$ means focusing in the $x$
direction, and $k<0$ means defocusing in the $x$ direction and focusing in
the $y$ direction), and for a sextupole $\rho =\infty $ and 
\begin{equation}
A_{s}=-\frac{h}{3}\left( x^{3}-xy^{2}\right) ,
\end{equation}
where $h$ is the sextupole strength. Combining these types of elements, we
set up a cell, which consists of the following sequence of elements: drift,
defocusing quadrupole superimposed with a sextupole, drift, bending dipole,
drift, defocusing quadrupole superimposed with a sextupole, and drift (Fig. 
\ref{fig:fodo}). The defocusing quadrupoles have the same strength $%
k=-0.0085,$ and the sextupole strength is $h=0.06$. The length of the drifts
are $1$ $m$, of the quadrupoles and the dipole $0.5$ $m$, and the dipole's
curvature radius $\rho =2.5$ $m$. Integration over this structure yielded a
Taylor model containing the true solution. Moreover, the final Taylor model
encloses the final conditions with a relative overestimation in the order of 
$10^{-10}$. Thus, the Taylor model approach does not only give verified
enclosures of the flow, but it also leads to only marginal overestimation in
the process. This result underscores the previously mentioned avoidance of
the wrapping effect by the Taylor model approach. As an example of this,
consider the Taylor model shown in Table \ref{table:tm}. By providing a
rigorous remainder bound, it describes a guaranteed enclosure of the final $%
x $ positions as a polynomial function of the initial conditions $x_{0}$, $%
a_{0}$, $y_{0}$, and $b_{0}$ at the end of the accelerator cell. For
increase readability, only terms of orders three or less have been listed
(the actual computations were done in order 17 and the reference polynomial
has a total of 1933 non-vanishing terms). Note that the initial regions are
all scaled to $[-1,1]$ internally (naturally, the polynomial coefficients
have been scaled accordingly).

\begin{figure}
\centering
\includegraphics[scale=0.3]{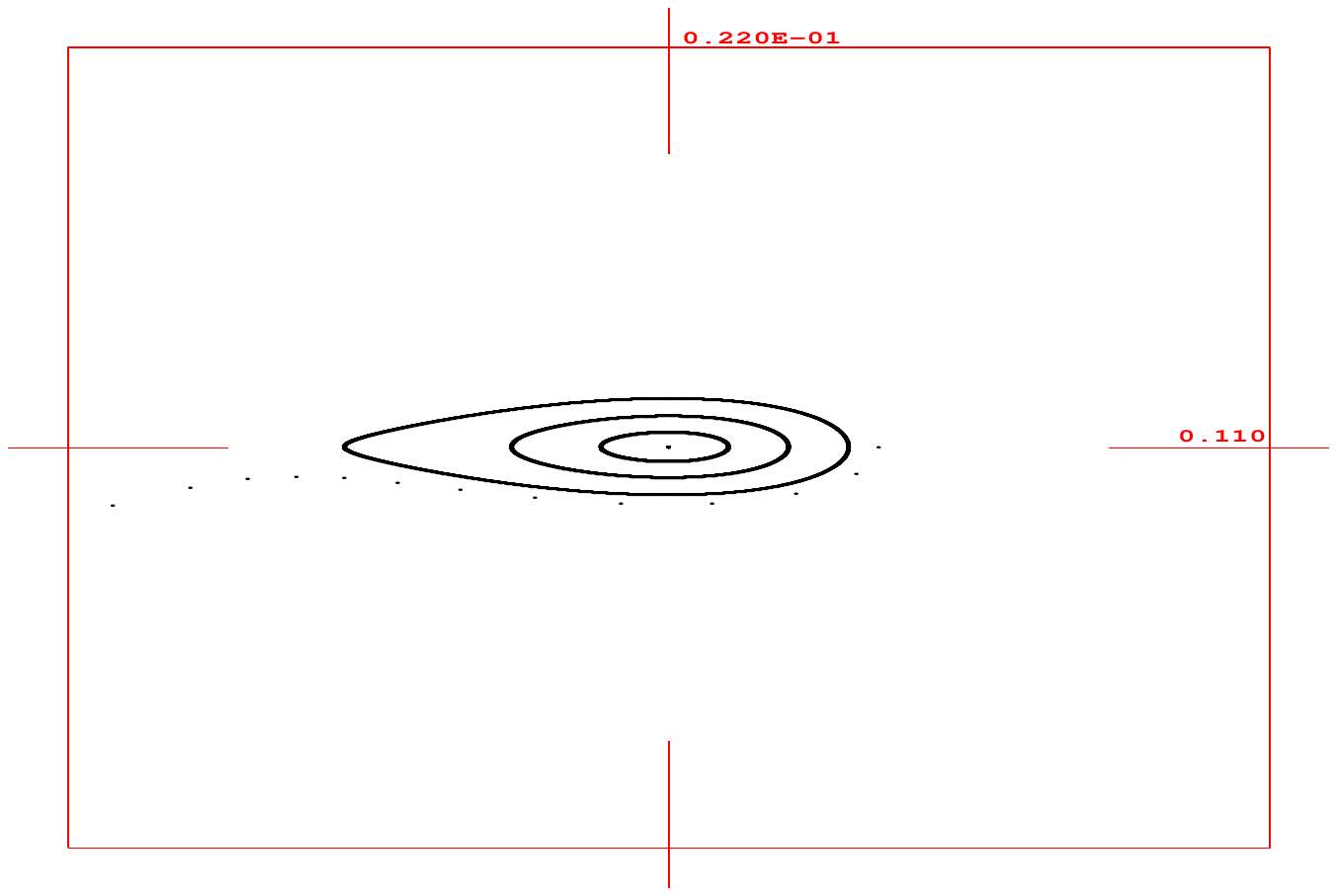}
\includegraphics[scale=0.3]{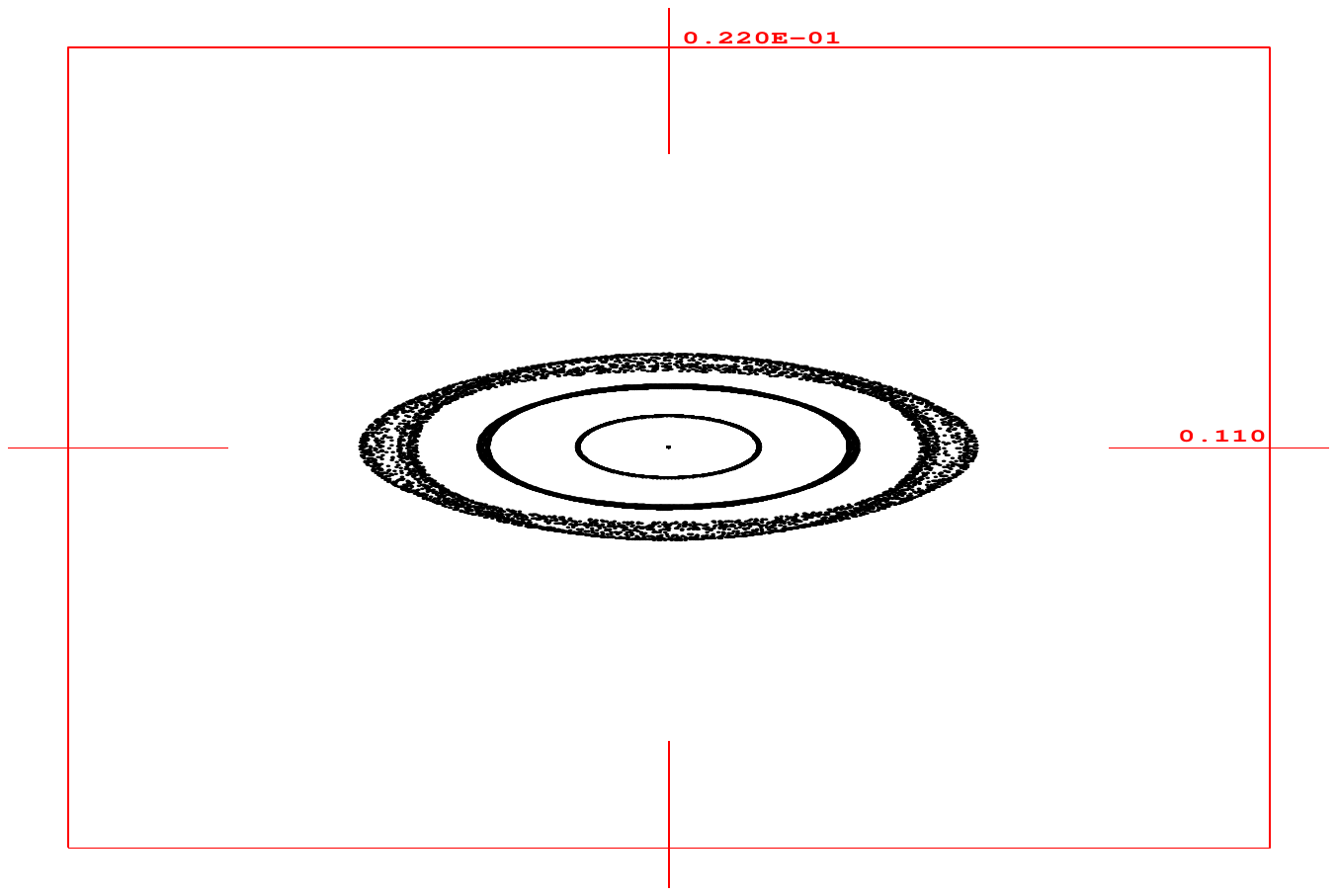}
\caption
{$(a)$ $(x,a)$, and $(b)$ $(y,b)$ tracking pictures of the one-turn map of 
the accelerator cell for particles launched along the $x$ and $y$ axes respectively. 
The box of guarantied invertibility of the generator associated with $S=0$ extends to at 
least $[-0.1,0.1]$ in every direction.}
\label{fig:celltrack}
\end{figure}

\begin{figure}
\centering
\includegraphics[scale=0.65]{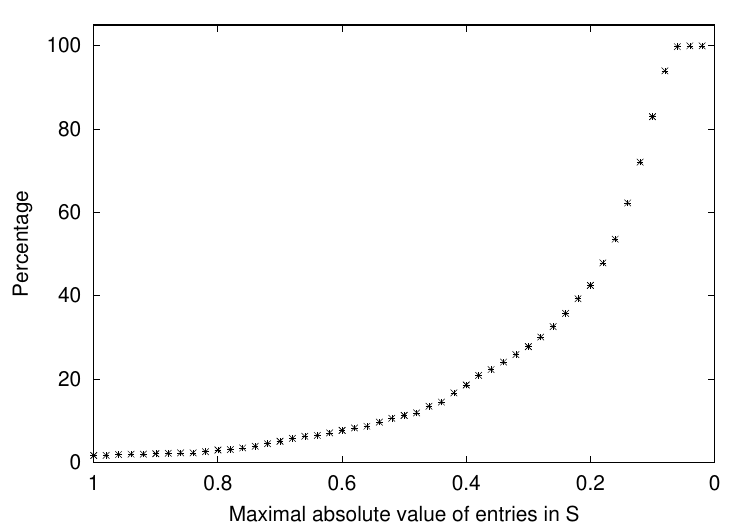}
\caption
{The size of the random symmetric matrices $S$ vs. the percentage of associated 
generators, which are defined at least up to $[-0.1,0.1]$ in every direction in phase 
space for the accelerator cell.}
\label{fig:size-of-s-fodo}
\end{figure}

As for the Fermi-Pasta-Ulam system, a plot was created, which shows the
percentage of a set of $1000$ random generators defined over a phase space
region of $\left[ -0.1,0.1\right] $ in every direction. The value of $0.1$
approximately corresponds to a typical accelerator magnet aperture, and
includes the dynamic aperture, which is estimated from tracking to be $%
\left( x,y\right) =\left( 0.03,0.045\right) $, as shown in Fig. \ref
{fig:celltrack}. Clearly, the trend is the same as for the Fermi-Pasta-Ulam
system: for smaller $S$ there is a larger chance that the associated
generating function is defined over the region of interest, and if $S$ is
close enough to zero, $100\%$ of the generators are defined (Fig. \ref
{fig:size-of-s-fodo}).

\begin{figure}
\centering
\includegraphics[scale=0.65]{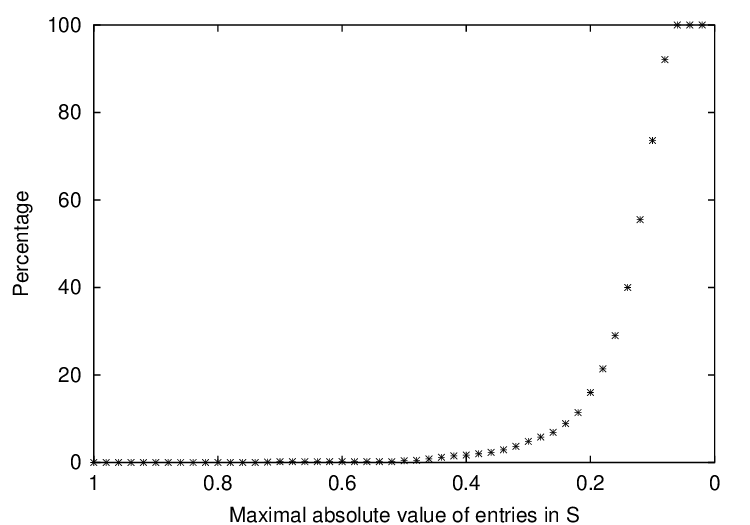}
\caption
{The size of the random symmetric matrices $S$ vs. the percentage of associated 
generators, which are defined at least up to $[-1,1]$ in every direction in phase space 
for the Fermi-Pasta-Ulam system.}
\label{fig:size-of-s}
\end{figure}

Finally, we note that the size of the domains of definition of the extended
generating functions is directly related to the problem of optimal
symplectic approximation of Hamiltonian flows \cite{mythesis}. The
conclusion reached in this paper are confirmed by a very general theoretical
argumentation based on Hofer's metric; in general the generating function
associated with $S=0$ gives the optimal approximation. This suggests that
the $S=0$ generator, as shown in Figs. \ref{fig:size-of-s} and \ref
{fig:size-of-s-fodo}, is the not the best because the Taylor model based
tools utilized to prove lower bounds for the domains of definition happen to
work best for small symmetric matrices $S$.

\bibliography{refs}

%TCIMACRO{
%\TeXButton{Fig-ex1}{\begin{figure}
%\centering
%\caption
%{Tracking picture of the cubic two dimensional symplectic map, and  the box of 
%guarantied invertibility of the generator associated with $S=0$.}
%\label{fig:2dpoly}
%\end{figure}%
%}}%
%BeginExpansion
%
%EndExpansion

%TCIMACRO{
%\TeXButton{Fig-ex2}{\begin{figure}
%\centering
%\caption
%{$(q_1,p_1)$ tracking picture of the four dimensional symplectic polynomial, 
%for two initial conditions launched along the $q_1$ axis close to the origin.}
%\label{fig:4dpoly-1}
%\end{figure}
%
%\begin{figure}
%\centering
%\caption
%{$(a)$ $(q_1,p_1)$, and $(b)$ $(q_2,p_2)$ tracking pictures of the four dimensional 
%symplectic polynomial for two particles (launched along the $q_1$ and $q_2$ axes 
%respectively) close to the dynamic aperture, and  the box of 
%guarantied invertibility of the generator associated with $S=0$.}
%\label{fig:4dpoly-2}
%\end{figure}%
%}}%
%BeginExpansion

%
%EndExpansion

%TCIMACRO{
%\TeXButton{Fig-ex3}{\begin{figure}
%\centering
%\caption
%{The Fermi-Pasta-Ulam system used in our example.}
%\label{fig:springs}
%\end{figure}
%
%\begin{figure}
%\centering
%\caption
%{$(a)$ $(q_1,p_1)$, and $(b)$ $(q_1,q_2)$ tracking pictures of the half period map of 
%the Fermi-Pasta-Ulam system for a particle launched along the $q_1$ axis. The box of 
%guarantied invertibility of the generator associated with $S=0$ extends to at least 
%$[-1,1]$ in every direction.}
%\label{fig:fpu}
%\end{figure}
%
%\begin{figure}
%\centering
%\caption
%{The size of the random symmetric matrices $S$ vs. the percentage of associated 
%generators, which are defined at least up to $[-1,1]$ in every direction in phase space 
%for the Fermi-Pasta-Ulam system.}
%\label{fig:size-of-s}
%\end{figure}%
%}}%
%BeginExpansion

%
%EndExpansion

%TCIMACRO{
%\TeXButton{Fig-ex4}{\begin{figure}
%\centering
%\caption
%{The accelerator cell used in the example of Subsection \ref{ssect:fodo}. The drifts are 
%$1 m$ long, the quadrupoles and the dipole are $0.5 m$ long, and the dipole's radius 
%of curvature is $2.5 m$. The first quadrupoles are defocusing in the $x$ direction, 
%with strength $k=-0.0085$, and the reference particle is a $1 MeV$ proton.}
%\label{fig:fodo}
%\end{figure}
%
%\begin{figure}
%\centering
%\caption
%{$(a)$ $(x,a)$, and $(b)$ $(y,b)$ tracking pictures of the one-turn map of 
%the accelerator cell for particles launched along the $x$ and $y$ axes respectively. 
%The box of guarantied invertibility of the generator associated with $S=0$ extends to at 
%least $[-0.1,0.1]$ in every direction.}
%\label{fig:celltrack}
%\end{figure}
%
%\begin{figure}
%\centering
%\caption
%{The size of the random symmetric matrices $S$ vs. the percentage of associated 
%generators, which are defined at least up to $[-0.1,0.1]$ in every direction in phase 
%space for the accelerator cell.}
%\label{fig:size-of-s-fodo}
%\end{figure}%
%}}%
%BeginExpansion

%
%EndExpansion

%TCIMACRO{
%\TeXButton{B}{\begin{table}
%\centering%
%}}%
%BeginExpansion
\begin{table}
\caption
{Lower bounds for the domains of definition of several generating function types for 
the two dimensional cubic symplectic polynomial. The matrices $S_1$ through $S_{10}$ 
are randomly generated with entries in $[-10,10]$.} 
\centering%
%
%EndExpansion
\vspace{1em}
\begin{tabular}{cc}
$\left[ S\right] $ & $\left( q_{1};p_{1}\right) $ lower bounds for the
domain of the generator \\ \hline
$0$ & $\left( [-0.333,0.333];[-0.333,0.333]\right) $ \\ 
$I$ & $\left( [-0.310,0.310];[-0.310,0.310]\right) $ \\ 
$S_{1}$ & $\left( [-0.301,0.301];[-0.301,0.301]\right) $ \\ 
$S_{2}$ & $\left( [-0.261,0.261];[-0.261,0.261]\right) $ \\ 
$S_{3}$ & $\left( [-0.258,0.258];[-0.258,0.258]\right) $ \\ 
$S_{4}$ & $\left( [-0.360,0.360];[-0.360,0.360]\right) $ \\ 
$S_{5}$ & $\left( [-0.377,0.377];[-0.377,0.377]\right) $ \\ 
$S_{6}$ & $\left( [-0.109,0.109];[-0.109,0.109]\right) $ \\ 
$S_{7}$ & $\left( [-0.111,0.111];[-0.111,0.111]\right) $ \\ 
$S_{8}$ & $\left( [-0.105,0.105];[-0.105,0.105]\right) $ \\ 
$S_{9}$ & $\left( [-0.090,0.090];[-0.090,0.090]\right) $ \\ 
$S_{10}$ & $\left( [-0.042,0.042];[-0.042,0.042]\right) $%
\end{tabular}
%TCIMACRO{
%\TeXButton{E}{\caption
%{Lower bounds for the domains of definition of several generating function types for 
%the two dimensional cubic symplectic polynomial. The matrices $S_1$ through $S_{10}$ 
%are randomly generated with entries in $[-10,10]$.}
%\label{table:2dpoly}
%\end{table}%
%}}%
%BeginExpansion

\label{table:2dpoly}
\end{table}%
%
%EndExpansion

%TCIMACRO{
%\TeXButton{B}{\begin{table}
%\centering%
%}}%
%BeginExpansion
\begin{table}
\caption
{Lower bounds for the domains of definition of several generating function types for 
the four dimensional symplectic polynomial. The matrices $S_1$ through $S_{10}$ 
are randomly generated with entries in $[-10,10]$.}
\centering%
\vspace{1em}
%
%EndExpansion
\begin{tabular}{cc}
$\left[ S\right] $ & $\left( q_{1};p_{1};q_{2};p_{2}\right) $ lower bounds
for the domain of the generator \\ \hline
$0$ & $\left(
[-0.02220,0.03404];[-0.01924,0.01924];[-0.02664,0.03552];[-0.05328,0.05328]%
\right) $ \\ 
$I$ & $\left(
[-0.03742,0.05738];[-0.03243,0.03243];[-0.04491,0.05988];[-0.08982,0.08982]%
\right) $ \\ 
$S_{1}$ & $\left(
[-0.01732,0.02656];[-0.01501,0.01501];[-0.02079,0.02772];[-0.04158,0.04158]%
\right) $ \\ 
$S_{2}$ & $\left(
[-0.01610,0.02474];[-0.01398,0.01398];[-0.01936,0.02582];[-0.03873,0.03873]%
\right) $ \\ 
$S_{3}$ & $\left(
[-0.02112,0.03238];[-0.01830,0.01830];[-0.02534,0.03379];[-0.05068,0.05068]%
\right) $ \\ 
$S_{4}$ & $\left(
[-0.04900,0.07514];[-0.04247,0.04247];[-0.05880,0.07840];[-0.11761,0.11761]%
\right) $ \\ 
$S_{5}$ & $\left(
[-0.00967,0.01483];[-0.00838,0.00838];[-0.01161,0.01548];[-0.02322,0.02322]%
\right) $ \\ 
$S_{6}$ & $\left(
[-0.01630,0.02500];[-0.01413,0.01413];[-0.01956,0.02608];[-0.03913,0.03913]%
\right) $ \\ 
$S_{7}$ & $\left(
[-0.01489,0.02283];[-0.01290,0.01290];[-0.01787,0.02383];[-0.03574,0.03574]%
\right) $ \\ 
$S_{8}$ & $\left(
[-0.01186,0.01819];[-0.01028,0.01028];[-0.01423,0.01898];[-0.02847,0.02847]%
\right) $ \\ 
$S_{9}$ & $\left(
[-0.02064,0.03164];[-0.01788,0.01788];[-0.02476,0.03302];[-0.04953,0.04953]%
\right) $ \\ 
$S_{10}$ & $\left(
[-0.06874,0.10540];[-0.05957,0.05957];[-0.08249,0.10999];[-0.16498,0.16498]%
\right) $%
\end{tabular}
%TCIMACRO{
%\TeXButton{E}{\caption
%{Lower bounds for the domains of definition of several generating function types for 
%the four dimensional symplectic polynomial. The matrices $S_1$ through $S_{10}$ 
%are randomly generated with entries in $[-10,10]$.}
%\label{table:4dpoly}
%\end{table}%
%}}%
%BeginExpansion

\label{table:4dpoly}
\end{table}%
%
%EndExpansion

%TCIMACRO{
%\TeXButton{B}{\begin{table}
%\centering%
%}}%
%BeginExpansion
\begin{table}
\caption
{Taylor model describing final $x$ positions as a function of the initial conditions 
$x_0,a_0,y_0,b_0$ for the accelerator cell of subsection \ref{ssect:fodo}.}
\centering%
\vspace{1em}
\begin{BVerbatim}
  RDA VARIABLE:   NO=  17, NV=     4
     I  COEFFICIENT            ORDER EXPONENTS 
     1  0.1398389113940111       1   1 0  0 0
     2  0.1038317686361456       1   0 1  0 0
     3  -.2447944264979660E-01   2   2 0  0 0
     4  -.1183394850192213E-01   2   1 1  0 0
     5  -.2119694344941219E-02   2   0 2  0 0
     6  0.2361409770673162E-01   2   0 0  2 0
     7  0.1212766097410308E-01   2   0 0  1 1
     8  0.2185093364668458E-02   2   0 0  0 2
     9  0.9944830763540945E-03   3   3 0  0 0
    10  0.8715481705011164E-03   3   2 1  0 0
    11  0.1933878768639235E-03   3   1 2  0 0
    12  0.3115908412345901E-04   3   0 3  0 0
    13  0.1238155786756443E-02   3   1 0  2 0
    14  -.1757616408645450E-03   3   0 1  2 0
    15  0.1233139677764291E-02   3   1 0  1 1
    16  0.8414517599464383E-04   3   0 1  1 1
    17  0.1954684434876060E-03   3   1 0  0 2
    18  0.4920452675843482E-04   3   0 1  0 2
  ----------------------------------------------
  VAR    REFERENCE POINT                  DOMAIN INTERVAL
    1  0.000000000000000      [-1.00000000 , 1.00000000 ]
    2  0.000000000000000      [-1.00000000 , 1.00000000 ]
    3  0.000000000000000      [-1.00000000 , 1.00000000 ]
    4  0.000000000000000      [-1.00000000 , 1.00000000 ]
  ORDER                    BOUND INTERVAL
    1     [-.2436706800301567     ,0.2436706800301567     ]
    2     [-.5056074647076302E-001,0.4976080054742529E-001]
    3     [-.5066453459468558E-002,0.5066453459468558E-002]
    R     [-.2157121190249145E-012,0.2178948979195422E-012]
  ********************************************************* 
\end{BVerbatim}
%
%EndExpansion
%TCIMACRO{
%\TeXButton{E}{\caption
%{Taylor model describing final $x$ positions as a function of the initial conditions 
%$x_0,a_0,y_0,b_0$ for the accelerator cell of subsection \ref{ssect:fodo}.}
%\label{table:tm}
%\end{table}%
%}}%
%BeginExpansion
\label{table:tm}
\end{table}%
%
%EndExpansion

\end{document}